\documentstyle[times,elsart12]{elsart}

\input amssym.def
\input amssym
\input xypic

\newtheorem{Theorem}{Theorem}[section]
\newtheorem{Lemma}[Theorem]{Lemma}

\newtheorem{Proposition}[Theorem]{Proposition}
\newtheorem{Remark}[Theorem]{Remark}

\newtheorem{Conjecture}[Theorem]{Conjecture}
\newtheorem{Definition}[Theorem]{Definition}
\newtheorem{Question}[Theorem]{Question}

\def\tratto{\mbox{\rule{2mm}{.2mm}$\;\!$}}

\newfont{\bella}{eusm10 scaled\magstep1}
\def\F{\mbox{\bella \symbol{'106}}}

\def\E{\mbox{\bella \symbol{'105}}}
\def\Z{\mbox{\bella \symbol{'132}}}
\def\B{\mbox{\bella \symbol{'102}}}

\def\N{\mbox{\bella \symbol{'116}}}
\def\Ko{\mbox{\bella \symbol{'113}}}

\newfont{\Bella}{eusm10 scaled\magstephalf}
\def\KKo{\mbox{\small{\Bella \symbol{'113}}}}

\newfont{\frank}{eufm10 scaled\magstep1}
\def\p{\mbox{\frank p}}
\def\q{\mbox{\frank q}}
\def\m{\mbox{\frank m}}

\newfont{\frankk}{eufm10 scaled\magstephalf}
\def\pp{\mbox{\frankk p}}
\def\qq{\mbox{\frankk q}}
\def\mm{\mbox{\frankk m}}

\begin{document}
\begin{frontmatter}
\title{\Large\bf Reduction Number of Links of
Irreducible Varieties\thanksref{CNR}}
\thanks[CNR]{Both authors gratefully acknowledge partial support
from the {\em Consiglio Nazionale delle Ricerche} $($Italy$)$, under CNR
grant 203.01.63.}
\author{Alberto Corso\thanksref{RU}} 
\qquad \author{Claudia Polini\thanksref{MSU}}
\address{Department of Mathematics, Rutgers University,
New Brunswick, NJ 08903 $($USA$)$}
\thanks[RU]{Corresponding author. E-mail address: corso@math.rutgers.edu.}
\thanks[MSU]{The current address of the author is:
Department of Mathematics - Michigan State University - East Lansing,
Michigan 48824. E-mail address: polini@math.msu.edu.}
\date{August 16, 1995}

\begin{abstract}
The reductions of an ideal $I$ give a natural pathway to the
properties of $I$, with the advantage of having fewer generators.
In this paper we primarily focus on a conjecture about the reduction
exponent of links of a broad class of primary ideals. The existence of
an algebra structure on the Koszul and Eagon--Northcott resolutions is
the main tool for detailing the known cases of the conjecture.
In the last section we relate the conjecture to a formula involving
the length of the first Koszul homology modules of these ideals.

\

\par\noindent {\it 1991 Math. Subj. Class.:} Primary 13H10; Secondary 13C40,
13D40, 13D45, 13H15
\end{abstract}
%
%
\end{frontmatter}

\section{Introduction}
Blowing up a variety along a subvariety is a very important kind of
transformation in Algebraic Geometry, especially in the process of
resolution of singularities. Algebraically, this is encoded in the
following diagram
\begin{equation}
{ \diagram
\ \, \mbox{Proj}(\mbox{gr}_I(R))\ \, \rto|<\hole|<<\ahook \dto|>>\tip &
\ \, \mbox{Proj}(R[It]) \dto^{\varphi} \ \, \\
\qquad V(I) \qquad \rto|<\hole|<<\ahook & \quad \mbox{Spec}(R) \quad
\enddiagram }\label{proj}
\end{equation}
where $R$ is a commutative Noetherian local ring and $I$ is one of its
ideals.
The existence of the morphism $\varphi$ in $($\ref{proj}$)$, which is an
isomorphism outside $\varphi^{-1}(V(I))=\mbox{Proj}(\mbox{gr}_I(R))$,
justifies the attention paid to find conditions for the {\em Rees
algebra} and the {\em associated graded ring} of $I$, collectively referred
to as {\em blowup algebras},
\[
R[It]=\bigoplus_{i=0}^{\infty}  I^it^i\qquad \mbox{and } \qquad
\mbox{gr}_I(R)=\bigoplus_{i=0}^{\infty} I^i/I^{i+1},
\]
to be Cohen--Macaulay.

In what follows, it is described how methods from Linkage
Theory $($see \cite{PS}$)$ are providing advances on the
Cohen--Macaulayness issue.
Recall that two ideals $I$ and $L$ of height $g$ in a Cohen--Macaulay local
ring $R$ are said to be {\em directly linked} if there exists a
regular sequence ${\bf z}= z_1, \ldots, z_g\subseteq
I\cap L$ such that $I= ({\bf z})\colon L$ and $L= ({\bf z})\colon I$.
Furthermore, they are said to be {\em geometrically linked} if they
are unmixed ideals, without common components
and $I\cap L={\bf z}$.

Also, a {\em reduction} of an ideal $I$ $($see \cite{NR}$)$ is an ideal
$J\subseteq I$ such that, for some non-negative integer $r$, the
equality $I^{r+1}=JI^r$ holds. The smallest such integer is the
reduction number $r_J(I)$ of $I$ relative to $J$. If the residue field
of $R$ is infinite, then minimal reductions
always exist and their minimal number of generators does not depend on
the minimal reduction. This number is called {\em analytic spread of $I$} and
is always bigger than or equal to the height of the ideal $I$; in case
of equality the ideal $I$ is said to be {\em equimultiple}.

The specific question motivating this paper is to find conditions for
a link of an irreducible variety to have Cohen--Macaulay blowup algebras.
Since irreducible varieties correspond to the locus of primary ideals,
the central object of this study consists of ideals of the form
$I=J \colon P$, where $J=({\bf z})$ is an ideal generated by a regular
sequence of length $g$ inside $P$ and $P$ is a primary ideal of height
$g$. Of course, one cannot expect a positive answer to the previous
question if $P$ is any primary ideal. The first issue is then to
single out families of primary ideals for which good results are
expected. The first two natural choices are: $(${\em a}$)$ links of prime
ideals; $(${\em b}$)$ links of symbolic powers of prime ideals.

Two previous papers $($see \cite{CP,CPV}$)$ were devoted to
the study of the first family of ideals and the key point in order to have a
good understanding of that family was to show that links of prime
ideals are equimultiple with reduction number one. Indeed, the
counterpart of a low reduction exponent is, in general, a very high
depth of the blowup algebras associated to the ideals $($see
\cite[Corollary 2.4]{CP} and \cite[Theorem 3.1, Corollary
3.2]{CPV}$)$.  In its broadest formulation one has

\begin{Theorem}\label{B}
Let $(R, \m)$ be a Cohen--Macaulay local ring and let $J=(z_1,
\ldots, z_g)$ be an ideal generated by a regular sequence inside a
prime ideal $\p$ of height $g$. If we set $I=J \colon \p$
then $I^2=JI$ if one of the following two conditions holds

\protect\noindent\parbox[t]{2pc}{$(${\rm L}${}_1)$}
\parbox[t]{30.7pc}{$R_{\pp}$ is not a regular local ring;}

\protect\noindent\parbox[t]{2pc}{$(${\rm L}${}_2)$}
\parbox[t]{30.7pc}{$R_{\pp}$ is a regular local ring of
dimension at least $2$ and two of the $z_i$'s lie in $\p^{(2)}$.}
\end{Theorem}

\noindent Theorem~\ref{B} is sharp in the sense that if condition $(${\rm
L}${}_2)$ is not satisfied then $I$ is generically a complete
intersection, hence it does not admit any proper reduction.
A straightforward generalization of Theorem~\ref{B} reads as follows

\begin{Theorem}
Let $(R, \m)$ be a Cohen--Macaulay local ring and let $P$ be
an unmixed radical ideal of $R$ of codimension $g$, such that
$R_{\pp}$ is not a regular local ring for any prime $\p$ minimal
over $P$. Then every link $I$ of $P$ is equimultiple with reduction
number one.
\end{Theorem}

The focus of this paper is on the second family of ideals,
called {\em $k$-iterated links of prime ideals} $($see
Definition~\ref{Definition2.3.2} below$)$; the main conjecture, whose
formulation is inspired by the one of Theorem~\ref{B}, is about
their reduction number and it can be stated as follows

\begin{Conjecture}\label{conj0}
Let $(R, \m)$ be a Cohen--Macaulay local ring and $\p$ a
prime ideal of height $g\geq 2$. Let $J=(z_1, \ldots, z_g)$ be an
ideal generated by a regular sequence contained in $\p^{(s)}$, where $s$
is a positive integer $\geq 2$. For any $k=1, \ldots, s$ set
$I_k=J\colon \p^{(k)}$. Then
\[
I^2_k=JI_k
\]
if one of the following two conditions holds

\protect\noindent\parbox[t]{2.8pc}{$(${\rm IL}${}_1)$}
\parbox[t]{29.9pc}{$R_{\pp}$ is not a regular local ring;}

\protect\noindent\parbox[t]{2.8pc}{$(${\rm IL}${}_2)$}
\parbox[t]{29.9pc}{$R_{\pp}$  is a regular local ring and two of the
$z_i$'s lie in $\p^{(s+1)}$.}

%
\end{Conjecture}

\noindent Note, however, that in the main result of this paper
$(${\rm IL}${}_2)$ is replaced by the weaker assumption

\protect\noindent\parbox[t]{2.8pc}{$(${\rm IL}${}_2^{\ast})$}
\parbox[t]{29.9pc}{\em $R_{\pp}$  is a regular local
ring and either its dimension is at least $3$ or $J\subseteq
\p^{(s+1)}$ whenever its dimension is $2$.}

\medskip

We now describe the results of this manuscript.
The main result of Section $2$ is Theorem~\ref{Theorem2.3.4}, which
establishes the above conjecture when: $(${\em a}$)$
condition $(${\rm IL}${}_2^{\ast})$ holds; $(${\em b}$)$ the {\em initial
forms} $($see \cite[page 179]{BH}$)$ of the
generators of $J$ form a regular sequence on the associated graded
ring of $\p$. The approach is ideal-theoretic as in the case of
\cite{CP}. The first step is the reduction to the case of
$\m$-primary ideals; then, standard mapping cone arguments lead to a
detailed analysis of the structure of the ideals $I_k$. The main tool
for this analysis is the existence of an algebra
structure on the Koszul and the Eagon--Northcott resolutions; this
accounts for hypothesis $(${\em a}$)$. Hypothesis $(${\em
b}$)$ is, on the other hand, necessary in order to use a well-established
criterion due to Valabrega--Valla $($see \cite[Corollary
2.7]{ValabregaValla}$)$.

In Section $3$, instead, $R_{\pp}$ is only supposed to be
Gorenstein. The approach, here, is more homological rather than
ideal-theoretic. More precisely, we relate the length of $I^2_k/JI_k$ to the
length of the kernel $\delta(I_k)$ of the natural surjection from
$S_2(I_k) $ onto $I_k^2$. Furthermore, we conjecture that the sum of
these two lengths equals ${n_k+1 \choose 2}$, where
$n_k=\lambda(\m^{k-1}/\m^k)$. This fact is already known in the case
$k=1$ $($see \cite{CPV}$)$; we now show that it holds in the case
$k=2$ as well. The proof is based on a formula for the length of the
first Koszul homology module of $I_2$ $($see Theorem~\ref{H_1}$)$.

As a general rule, we will provide the basic definitions and the
references for all the results that are used; the reader
desiring more details should then consult the invaluable books of
W. Bruns and J. Herzog \cite{BH}, H. Matsumura \cite{Mat}, and the
excellent monograph of W. V. Vasconcelos \cite{WVV}.

\section{Iterated links of prime ideals in regular local rings}
A useful method to test whether a link of an ideal is equimultiple
with reduction number one is stated in the next lemma.

\begin{Lemma}\label{Lemma2.3.1}
Let $(R, \m)$ be a local ring and let $I$ be one of its ideals.
If $I$ is an equimultiple ideal of height $g$ with reduction number
$1$, i.e., $I^2=JI$ where $J$ is an ideal generated by a regular
sequence $z_1, \ldots, z_g$, then
\begin{equation}
I\subseteq J\colon I \qquad \mbox{and} \qquad
(J\colon I)J=(J\colon I)I. \label{eq1}
\end{equation}
Conversely, the conditions in $(\ref{eq1})$ are sufficient to
guarantee that $I$ has reduction number $1$ with respect to $J$
provided $I=J \colon (J \colon I)$.
\end{Lemma}
\begin{pf}
From $I^2=JI\subseteq J$ it follows that $I\subseteq J \colon I$, and
this takes care of the first condition in $($\ref{eq1}$)$---$J$ does
not need to be generated by a regular sequence.
Since $(J\colon I)J$ is clearly contained in $(J\colon I)I$, it is enough
to show that any $z\in (J\colon I)I$ can be written as an element in
$(J\colon I)J$. But $(J\colon I)I \subseteq J$, so that there exist
$\alpha_1, \ldots, \alpha_g$ such that
\begin{equation}\label{alpha}
z=\sum_{i=1}^g \alpha_i z_i.
\end{equation}
The proof will be complete provided it is shown that $\alpha_i$ is in $J
\colon I$ for all $i$ in the range $1, \ldots, g$. Pick $a\in I$ and
consider the element $za$ in $(J\colon I)I^2=(J\colon I)JI\subseteq
J^2$; Hence
\begin{equation}
za=\sum_{i=1}^g \beta_i z_i, \label{beta}
\end{equation}
for some $\beta_1, \ldots, \beta_g$ in $J$. A comparison between
equation $(\ref{beta})$ and $a$ times equation $(\ref{alpha})$ yields
the following identity
\[
\sum_{i=1}^g (\alpha_i a -\beta_i)z_i=0,
\]
which says that $(\alpha_i a -\beta_i)z_i \equiv 0 \ \mbox{mod} \ (z_1,
\ldots, \widehat{z_i}, \ldots, z_g)$, i.e., $\alpha_i a \in J$ for all
$i$ as claimed.

The proof of the converse is similar and can be found, essentially, in
\cite[Theorem 2.2]{CP}.  \qed
\end{pf}

\begin{Definition}\label{Definition2.3.2}
Let $(R, \m)$ be a Cohen--Macaulay local ring and $\p$ a prime
ideal of height $g$. Let $J$ be an ideal generated by a regular
sequence contained in $\p^{(s)}$, where $s$ is a positive integer.
Given any integer $k=1, \dots, s$, define $I_k=J\colon
\p^{(k)}$ and call it {\em $k$-th iterated link} of $\p$, as
$I_k=I_{k-1}\colon \p$.
\end{Definition}

\begin{Lemma}\label{Lemma2.3.3}
In the situation described in {\rm Definition~\ref{Definition2.3.2}}
one has that

\protect\noindent\parbox[t]{2pc}{$($\mbox{a}$)$}\parbox[t]{30.7pc}{if
$J\subseteq \p^s$ then for $k=1, \ldots,
s$ $I_k=J\colon \p^k=J\colon \p^{(k)}$;}

\parbox[t]{2pc}{$($\mbox{b}$)$}\parbox[t]{30.7pc}{if $R_{\pp}$
is a Gorenstein local ring, then for $k=1, \ldots, s$
$I_k$ and $\p^{(k)}$ are directly linked via $J$.}
%
%
\end{Lemma}
\begin{pf}
$(${\em a}$)$ Set $H_k=J\colon \p^k$. Then $\p^k I_k \subseteq \p^{(k)} I_k
\subseteq J$ says that $J\colon \p^{(k)}=I_k \subseteq
H_k=J\colon \p^k$. Hence, it suffices to show the
equality at the associated primes of the ideal $I_k$, which are
contained in the ones of $J$; hence they all have height $g$.
The conclusion now follows as $\p^{(k)}=\p^kR_{\pp}
\cap R$. This takes care of $(${\em a}$)$.

$(${\em b}$)$ This is a basic property of Linkage Theory $($see
\cite{PS}$)$. \qed
\end{pf}

Next, we formulate the main result of this paper; the rest of this
section will be entirely devoted to its proof, which will be
achieved in two steps. First, using a well-established method $($see
\cite{CP,CPV}$)$, the problem is reduced to the maximal case, i.e., to
ideals of the form $I_k=J \colon \m^k$.
Subsequently, using Remark~\ref{Remark2.3.8} together with a mapping
cone construction which produces a free resolution of $I_k$, one is
able to outline the set-theoretic structure of the ideal $I_k$. More
precisely, in Proposition~\ref{Proposition2.3.10} it is shown that
$I_k$ is made up of $J$ and extra ${d+k-2\choose d-1}$ generators
that lie in a sufficiently large power of $\m$ $($here $d$ denotes the
dimension of the ring $R)$. The notions of order
and initial form of an element are central in the following result so
they are formalized next.
The {\em order} $o(r)$ of any element $r$ in $R$ is defined as
the unique positive integer $q$ such that
$r\in \m^q$ but $r\not\in \m^{q+1}$. If $o(r)=q$, then the {\em
initial form} of $r$ is the residue class of $r$ in $\m^q/\m^{q+1}$.

\begin{Theorem}\label{Theorem2.3.4}
Let $(R, \m)$ be a Cohen--Macaulay local ring and $\p$
a prime ideal of height $g\geq 2$ such that $R_{\pp}$ is a regular local
ring. Let $J=(z_1, \ldots, z_g)$ be an ideal generated by a regular
sequence contained in $\p^{(s)}$, where $s$ is a positive integer $\geq
2$, and assume that the initial forms of the $z_i$'s form a
regular sequence on $\mbox{gr}_{\pp}(R)$. Then for any $k=1, \ldots, s$
the ideal $I_k=J \colon \p^{(k)}$
satisfies
\[
I^2_k=JI_k
\]
if the height of $\p$ is at least $3$ or if $J\subseteq \p^{(s+1)}$
whenever the height of $\p$ is $2$.
\end{Theorem}

%
%

\begin{Remark}[Reduction to the maximal case]\label{Remark2.3.6}
In {\rm Theorem~\ref{Theorem2.3.4}} $\p$ may be assumed to be the
maximal ideal of the ring.
{\rm
\begin{pf} As seen in \cite{CP,CPV}, it is enough to establish the
equality $I^2_k=JI_k$ at the associated primes of $JI_k$. From the
exact sequence
\[
\diagram
0 \rto & J/JI_k  \rto & R/JI_k \rto & R/J \rto & 0
\enddiagram
\]
and the fact that $J/JI_k=J/J^2\otimes R/I_k\simeq (R/I_k)^g$, it
follows that the associated primes of $JI_k$ are contained in those
of $R/J$ and $R/I_k$. On the other hand, $I_k$ and $\p^{(k)}$
are directly linked via $J$ so that $\mbox{Ass}(R/I_k)\subseteq R/J$.
This implies that the associated primes of $JI_k$ must have height
$g$. Therefore for $\q \in \mbox{Ass}(R/J)$ either $\q=\p$ or
$\q\not\supseteq \p$.
In the second case $(I_k)_{\qq}=J_{\qq}$ and the proof is complete;
the first case instead is the case of a link with a power of the
maximal ideal. Clearly, all the other conditions localize. \qed
\end{pf} }
\end{Remark}

\begin{Remark}\label{Remark2.3.7}
{\rm A resolution $\F=\{\F_i\}$ with differentials
$\partial=\{ d_i\}$ of a cyclic module $R/H$ is called an {\em algebra
resolution} if there is a graded associative multiplication $($called
a {\em {\rm DG}-algebra structure}$)$ $\diagram \F\otimes_R
\F \rto & \F \enddiagram$ lifting the usual product on $R/H$
and satisfying
\begin{enumerate}
\item $x_ix_j=(-1)^{ij}x_jx_i$ for $x_i\in \F_i$ and $x_j\in \F_j$;

\item $d_{i+j}(x_ix_j)=(d_ix_i)x_j+(-1)^ix_i(d_jx_j)$ for $x_i\in \F_i$
and $x_j\in \F_j$.
\end{enumerate}
If also $\partial\, \F\subset {{\m}}\F$ then $\F$ is called a
{\em minimal algebra resolution}.
The best known example is the Koszul resolution, which has
the structure of an exterior algebra on a free module. }
\end{Remark}

\begin{Remark}\label{Remark2.3.8}
{\rm
Set $d$ to be the dimension of the ring $R$.
Let $\Ko$ and $\E\N$ denote the Koszul resolution of $J=(z_1, \ldots,
z_d)$ and the Eagon--Northcott resolution of $\m^k$, respectively. The
natural inclusion $J \subseteq \m^k$ induces the following comparison
diagram
\begin{equation}\label{diagramma}{
\diagram
\E\N\, \colon\ & 0 \rto & \F_{k,d} \rto^{\delta_{d-1}}
\ddotted_{u_{k,d}}|<{\rotate\tip}
& \F_{k,d-1} \rto^{\delta_{d-2}}
\ddotted^{u_{k,d-1}}|<{\rotate\tip}&
\cdots\cdots  \rto^{\delta_1}  & \F_{k,1} \rto^{\delta_0}
\ddotted_{u_{k,1}}|<{\rotate\tip} &
\F_{k,0} \rto & 0  \\
\KKo\, \colon \ & 0 \rto & \Ko_d \rto_{d_d} & \Ko_{d-1}
\rto_{d_{d-1}} & \cdots \cdots
\rto_{d_2} & \Ko_1  \rto_{d_1} & \Ko_0 \rto
\uto_{u_{k,0}} & 0
\enddiagram }
\end{equation}
where
\[
\Ko_t=\bigwedge^t R^d \quad \ \mbox{for\ } \quad t\geq 0
\]
\[
\F_{k,0}=\bigwedge^k R^k
\]
\[
\F_{k,t}=D_{t-1}((R^k)^{\ast})\otimes_R \bigwedge^{k+t-1} R^{d+k-1}
\quad \ \mbox{for\ } \quad t\geq 1.
\]
Here, $D_{t-1}((R^k)^{\ast})$ denotes the degree $t-1$ component of
the {\em divided powers algebra} of the module $(R^k)^{\ast}$.

As the Koszul and the Eagon--Northcott resolutions both have an
associative, commutative, differential graded structure and the Koszul
resolution is a universal object, it is enough to define the map
$u_k=\{u_{k,t}\}$ on a basis of the Koszul resolution. To this end,
let $E_1, \ldots, E_d$ be a basis of $\Ko_1 \simeq R^d$ and $e_1,
\ldots, e_{d+k-1}$ be a basis of $R^{d+k-1}$. Then
\begin{eqnarray}\label{prodotto}
u_{k,d}(E_1\wedge \ldots \wedge E_d) & = &
u_{k,1}(E_1)\star \cdots \star u_{k,1}(E_d) \nonumber \\
& = & \prod\limits_{j=1}^d \sum\limits_{{\tiny \begin{array}{c} h=1 \\
\mid\! A_{h,k}\!\mid\, =k \end{array}}}^{{d+k-2\choose d-1}}
\alpha_{hj}(1\otimes e_{A_{h,k}}),
\end{eqnarray}
where $u_{k,1}$ is the map that rewrites the $z_i$'s in terms of the
basis elements of $\m^k$. Also, $A_{h,k}=(a_1, \ldots, a_k)$ is an
ordered $k$-upla with each $a_i$ an element of $\{1, \ldots, d+k-1\}$
and $e_{A_{h,k}}=e_{a_1} \wedge \ldots \wedge e_{a_k}$. }
\end{Remark}

%
%

\begin{Proposition}\label{Proposition2.3.10}
Let $(R, \m)$ be a regular local ring of dimension $d\geq 2$. Let
$J=(z_1, \ldots, z_d)$ be an ideal generated by a regular sequence in
$\m^s$, where $s$ is a positive integer $\geq 2$.
For any $k=1, \ldots, s$ set $I_k=J\colon \m^k$ and $q_i=o(z_i)$.
Then one can write $I_k=(J, J_k)$ with
\[
\mu(J_k)={d+k-2\choose d-1}
\qquad \mbox{and} \qquad  J_k\subseteq \m^{q-k+1},
\]
where $q=\sum\limits_{j=1}^d q_j -d$ if the dimension of $R$ is at
least $3$, or $q=q_1+q_2$ if the dimension of $R$ is $2$ and
$J\subseteq \m^{s+1}$.
\end{Proposition}
\begin{pf}
From $($\ref{diagramma}$)$, the fact that $I_k=J \colon \m^k$, and
\cite[Theorem 4.1.2]{WVV} it follows that the beginning of a free
resolution of $R/I_k$ is given by
\[
\diagram
\rto & (\Ko_{d-2}\oplus \F_{k,d-1})^{\ast}
\rto^{\partial_{d-1}^{\ast}} & (\Ko_{d-1}\oplus
\F_{k,d})^{\ast} \rto^(.6){\partial_d^{\ast}} & (\Ko_d)^{\ast}
\rto & R/I_k \rto & 0.
\enddiagram
\]
Since $\partial_d=(-d_d, -u_{k,d})^{{\rm T}}$, one has
\[
I_k=\partial_d^{\ast}((\Ko_{d-1}\oplus \F_{k,d})^{\ast})=(J,
-u_{k,d}(\Ko_d)^{{\rm T}}).
\]
Moreover, $\F_{k,d} \simeq R^{{d+k-2\choose d-1}}$ hence
$J_k=-u_{k,d}(\Ko_d)^{{\rm T}}$ has $\mu_k={d+k-2\choose d-1}$
generators, i.e.,
\begin{equation}\label{j.l}
I_k=(J, b_{k,1}, \ldots, b_{k,\mu_k}).
\end{equation}
Using $($\ref{prodotto}$)$, the fact that $\alpha_{hj}\in
\m^{q_j-k}$ $($or $\alpha_{hj}\in \m^{q_j+1-k}$ if $d=2)$,
and $d-1$ times the explicit formula for the product \,$\star$ given in
\cite[Section 6]{Hema}---which is valid regardless of the
characteristic of $R$---one concludes that the $b_{k,j}$'s belong to
the power $\m^e$, where
\[
e=\sum\limits_{j=1}^d q_j-dk+(d-1)(k-1)=\sum\limits_{j=1}^d q_j-d-k+1
\]
$($or $e=q_1+q_2-k+1$ if $d=2)$.
Thus, one has that $J_k \subseteq \m^{q-k+1}$, where
$q=\sum\limits_{j=1}^d q_j-d$ $($or $q=q_1+q_2$ if $d=2)$. \qed
\end{pf}

\begin{Remark}\label{Remark2.3.11}{\rm
Note that $q-k+1\geq q_d+(d-1)s-d-k+1=q_d+(d-1)(s-1)-k$. Thus
$q-k+1\geq q_d+s-k$ if $d\geq 3$, and $q-k+1=q_2+s-k+1> q_2+s-k$ if
$d=2$. }
\end{Remark}

\noindent We now have all the tools to complete the proof of
Theorem~\ref{Theorem2.3.4}.

\bigskip

\begin{pf*}{Proof of Theorem~\ref{Theorem2.3.4}.}
By a result of Valabrega--Valla $($see \cite[Corollary
2.7]{ValabregaValla}$)$, the elements of $J$ satisfy the following
equality for any $n\geq 1$
\[
\m^n\cap J=\sum_{i=1}^d \m^{n-q_i}z_i,
\]
where $q_i=o(z_i)$. As in a local ring regular sequences
permute, one may assume that $q_1\leq \cdots\leq q_d$. Hence
for any $n\geq 1$ it follows that
\begin{equation}\label{valabrega-valla}
\m^n\cap J\subseteq \m^{n-q_d}J.
\end{equation}
In particular, $($\ref{valabrega-valla}$)$ used with $n=q_d+s$
gives the inclusion
\begin{equation}\label{uva}
\m^{q_d+s}\cap J \subseteq \m^s J\subseteq \m^kJ.
\end{equation}
This is enough to prove the assertion. Indeed,
Proposition~\ref{Proposition2.3.10} and Remark~\ref{Remark2.3.11}
imply that $\m^kJ_k\subseteq \m^{q_d+s}$. On the other hand, from the
definition of $I_k$ it follows that $\m^kJ_k\subseteq \m^kI_k\subseteq
J$; hence $($\ref{uva}$)$ says that $\m^kJ_k\subseteq \m^k J$ and so
\[
\m^kI_k=\m^k(J, J_k)=\m^kJ+\m^kJ_k\subseteq \m^kJ.
\]
This completes the proof of the theorem as both the conditions
in Lemma~\ref{Lemma2.3.1} $($see $($\ref{eq1}$))$ are satisfied. \qed
\end{pf*}

\medskip

\begin{Question}{\rm
Let $I=J\colon H$, where $H$ is a perfect ideal with a minimal free
resolution having a DG-algebra structure $($see \cite{Matt} for a list
of such ideals$)$. Is it true that under some mild hypotheses $I^2=JI$?
Mild hypotheses on $J$ and $H$ are needed as the next example shows.

Let $R$ be a two dimensional regular local ring with maximal ideal
$\m=(x, y)$ and pick any positive integer $r\geq 2$. The $\m$-primary
ideals $J=(x^{r+1}, y^{r+1})$ and $H=(x, y^r)$ yield the following
commutative diagram
\[
\diagram
0 \rto & R \rrto^{-y^r \choose \ x} & & R^2 \rrto^{(x\ y^r)} & & R
\rto & R/H \rto & 0 \\
0 \rto & R \rrto_{-y^{r+1}\choose \ x^{r+1}}\uto^{\wedge^2\varphi}
& & R^2 \rrto_{(x^{r+1}\ y^{r+1})}\uto^{\varphi} & & R \rto\uto^{{\rm
id}} & R/J \rto\uto & 0
\enddiagram
\]
where $\varphi={x^r \ \ 0 \choose 0 \ \ \ y}$ is the matrix that writes
the generators of $J$ in terms of the ones of $H$, and
$\displaystyle{\wedge^2\varphi}=\det(\varphi)=x^ry$.
Both the deleted resolutions have a DG-algebra structure, as they are
Koszul resolutions.  However, the link $I=J\colon H=(x^{r+1}, y^{r+1},
x^ry)$ satisfies $r_J(I)=r \geq 2$. }
\end{Question}

\section{Iterated links of prime ideals in Gorenstein rings}
In this section we approach Conjecture~\ref{conj0} from a different
perspective, namely the one of \cite{CPV}.
In the case $k=1$, it turns out that $I_1=J\colon \m=(J, v_1,
\ldots, v_t)$, where $\overline{v}_1, \ldots, \overline{v}_t$ are the
generators of the socle of $R/J$ and $t$ is the {\em type} of $R$
$($see \cite[Definition 1.2.15]{BH}$)$. Moreover,
one also has the formula $($see \cite[Remark 2.7]{CP}$)$
\begin{equation}\label{jun}
\lambda(I_1^2/JI_1)+\lambda(\delta(I_1))={t+1\choose 2},
\end{equation}
where $\delta(I_1)$ is the kernel of the natural surjection from
$S_2(I_1)$ onto $I_1^2$; $(\ref{jun})$ is the key tool in the proof
of \cite[Theorem 2.1]{CPV}, since in the Gorenstein
case $t=1$ and the case $\delta(I_1)=0$ can not occur provided that
either condition $(${\rm L}${}_1)$ or $(${\rm L}${}_2)$ is satisfied.

Let $(R, \m)$ be a Gorenstein local ring of dimension $d$ and $J$ an
ideal of $R$ generated by a regular sequence $(z_1, \ldots, z_d)$
contained in $\m^s$, for some positive integer $s\geq 2$. Set
$I_k=J\colon \m^k$ for $k=1, \ldots, s$.
The goal of this section is to point out the analogue of
$($\ref{jun}$)$ for the $I_k$'s with $k\geq 2$. More precisely, we will
show in Proposition~\ref{Theorem2.4.8} that
\[
\lambda(I_k^2/JI_k)+\lambda(\delta(I_k))={n_k+1\choose 2}
\]
is equivalent to
\[
\lambda(H_1(I_k))=n_k\lambda(R/I_k)+{n_k+1\choose 2}-
\lambda(R/\m^k),
\]
where $n_k=\lambda(\m^{k-1}/\m^k)$.
We end the section by showing that the latter formula
holds in the case $k=2$ $($see Theorem~\ref{H_1}$)$. These are the
technical results of the manuscript and they are interesting
as they uncover an error in a formula that appears in
\cite{AA} and saying, in our terminology, that
\[
\lambda(H_1(I_k))\leq n_k \lambda(R/I_k).
\]

As, in what follows, $R$ need not be a regular local ring but
simply Gorenstein, one cannot use Proposition~\ref{Proposition2.3.10}
in order to describe the structure of $I_k$, since $R/\m^k$ does not
necessarily have a finite projective resolution. This is taken care of in
Proposition~\ref{lem1} below. Some additional assumptions will be
needed.

\begin{Remark}
{\rm
From now on, it will always be required that $\m^k\colon \m=\m^{k-1}$
for $k=3, \ldots,s$; this is always
satisfied if, e.g., $\mbox{depth}(\mbox{gr}_{\mm}(R))\geq 1$. Also,
note that $d\geq 2$ $($or, in the general case, $g\geq2)$ is needed as
the next example shows.

Let $R=k[X, Y, Z]/(Y^2-XZ, X^3-Z^2)=k[x, y, z]$, where $x$,
$y$, and $z$ denote the images of $X$, $Y$, and $Z$ modulo
$(Y^2-XZ, X^3-Z^2)$. It can be shown $($see \cite{Sally}$)$ that $R$
is a Gorenstein ring of dimension $1$ with $\mbox{gr}_{\mm}(R)$ Gorenstein
as well. Let $J=(x^4)$ be a regular sequence in $\m^4=(x, y, z)^4$ and
consider $I_3=J\colon \m^3$; a computation using the computer system
{\em Macaulay} shows that $I_3^2\neq JI_3$. }
\end{Remark}

\begin{Proposition}\label{lem1}
Let $(R, \m)$ be a Gorenstein local ring of dimension $d\geq 2$
and let $J=(z_1, \ldots, z_d)$ be an ideal generated by a regular
sequence contained in $\m^s$, where $s$ is a positive integer $\geq 2$.
For $k=1, \ldots, s$ set $I_k=J \colon \m^k$ and, for $k=3, \ldots,
s$, assume that $\m^k \colon \m = \m^{k-1}$. Then

\protect\noindent\parbox[t]{2pc}{$($\mbox{a}$)$}\parbox[t]{30.7pc}{$I_k/I_{k-1}$
is an $R/\m$-vector space of dimension $n_k$, where
$n_k=\lambda(\m^{k-1}/\m^k)$. Hence $I_k=(I_{k-1}, b_1,
\ldots, b_{n_k})$, with $\widehat{b}_1, \ldots,
\widehat{b}_{n_k}$ minimal generators of $\widehat{I}_k=I_k/I_{k-1}$;}

\protect\noindent\parbox[t]{2pc}{$($\mbox{b}$)$}\parbox[t]{30.7pc}{$I_k=(J,
b_1, \ldots, b_{n_k})$, with $\overline{b}_1, \ldots,
\overline{b}_{n_k}$ minimal generators of $\overline{I}_k=I_k/J$;}

\protect\noindent\parbox[t]{2pc}{$($\mbox{c}$)$}\parbox[t]{30.7pc}{$J$
is among the minimal generators of $I_k$ if and only if $\mu(I_k)=d+n_k$;}

\protect\noindent\parbox[t]{2pc}{$($\mbox{d}$)$}\parbox[t]{30.7pc}{$I_{k-1}=\m
I_k+J$.}
%
%
%
%
%
\end{Proposition}
\begin{pf}
One can write $I_k= I_{k-1} \colon \m$ so that $\m I_k\subseteq
I_{k-1}$. This shows that $I_k/I_{k-1}$ is an $R/\m$-vector space.
As $\mbox{Hom}_{R/J}(\tratto, R/J)$ is a length preserving functor $($see
\cite[Theorem 3.2.12]{BH}$)$, from the equality
$\mbox{Hom}_{R/J}(R/\m^k, R/J)=J\colon \m^k/J=I_k/J$ it follows that
\begin{equation}\label{11}
\lambda(I_k/J)=\lambda(\mbox{Hom}_{R/J}(R/\m^k, R/J))=\lambda(R/\m^k).
\end{equation}
Hence
\begin{equation}\label{1}
\lambda(I_k/I_{k-1})=\lambda(I_k/J)-\lambda(I_{k-1}/J)
=\lambda(\m^{k-1}/\m^k)=n_k.
\end{equation}
This completes the proof of $(${\em a}$)$.

$(${\em b}$)$ follows from the fact that $\mu(I_k/J)$ equals the type of
$R/\m^k$ $($see \cite[Proposition 3.3.11 $(${\em c},i$)$]{BH}$)$. More
precisely,
\begin{eqnarray*}
\mu(I_k/J) & = & \mbox{dim}(\mbox{Hom}_R(R/\m, R/\m^k)) \\
& = & \mbox{dim}(\m^k\colon \m/\m^k)=\lambda(\m^{k-1}/\m^k)=n_k.
\end{eqnarray*}
Note that in the last equality we used the hypothesis that $\m^{k-1}
\colon \m = \m^{k-1}$.

Finally, $(${\em c}$)$ and $(${\em d}$)$ readily follows from $(${\em
a}$)$ and $(${\em b}$)$. \qed
\end{pf}

\begin{Proposition}\label{Theorem2.4.8}
Let $(R, \m)$ be a Gorenstein local ring of dimension $d\geq 2$
and let $J=(z_1, \ldots, z_d)$ be an ideal generated by a regular
sequence contained in $\m^s$.
For $k=1, \ldots, s$ set $I_k=J \colon \m^k$ and, for $k=3, \ldots,
s$, assume that $\m^k \colon \m = \m^{k-1}$. Furthermore,
assume that the $z_i$'s are among the minimal generators of
$I_k$ and let $n_k=\lambda(\m^{k-1}/\m^k)$.
Then the following two results hold

\protect\noindent\parbox[t]{2pc}{$($\mbox{a}$)$}\parbox[t]{30.7pc}{the
length of the symmetric square $S_2(I_k/J)$ is given by
\begin{equation}\label{uv}
\lambda(S_2(I_k/J))=\lambda(I_k^2/JI_k)+\lambda(\delta(I_k))=
\rho ,
\end{equation}
where $\delta(I_k)$ is the kernel of the natural surjection from
$S_2(I_k)$ onto $I_k^2$;}

\protect\noindent\parbox[t]{2pc}{$($\mbox{b}$)$}\parbox[t]{30.7pc}{the
length of the first Koszul homology module $H_1(I_k)$ is given by
\begin{equation}\label{lh1}
\lambda(H_1(I_k))=n_k\lambda(R/I_k)+\rho-
\lambda(R/\m^k).
\end{equation}}
%
%
%

Moreover, $\rho={n_k+1 \choose 2}$ if and only if $S_2(I_k/J)$ is an
$R/\m$-vector space.
%
%
\end{Proposition}
\begin{pf}
As in the proof of \cite[Remark 2.7]{CP} one has the
short exact sequence
\[
\diagram
0 \rto & \delta(I_k) \rto & S_2(I_k/J) \rto & I_k^2/JI_k \rto & 0,
\enddiagram
\]
which leads to the equation
\begin{equation}\label{s2}
\lambda(S_2(I_k/J))=\lambda(I_k^2/JI_k)+\lambda(\delta(I_k)).
\end{equation}
If we set $\lambda(S_2(I_k/J))=\rho$, $(\ref{lh1})$ follows from
$(\ref{uv})$, \cite[Proposition 2.2]{CPV}, and $($\ref{11}$)$.

By Proposition~\ref{lem1} one has that $I_{k-1}=\m I_k+J$. Hence
the last assertion follows from
\[
S_2(I_k/I_{k-1}) = S_2(I_k/(\m I_k+J)) = S_2(I_k/J\otimes R/\m) \simeq
S_2(I_k/J)\otimes R/\m
\]
and the fact that $I_k/I_{k-1}$ is an $R/\m$-vector space of dimension
$n_k$. \qed
%
%
%
%
\end{pf}

\begin{Remark}
{\rm Note that the conjectured value for
$\rho=\lambda(S_2(I_k/J))$ does not depend on $I_k$ nor on $J$ but on
$\m$ alone. This follows from the fact that $I_k/J$ is isomorphic to
the canonical module of $R/\m^k$.}
\end{Remark}

\begin{Remark}\label{remh_1}
{\rm It also has to be pointed out that $($\ref{lh1}$)$ can be derived
in a more direct manner. It is well-known that $H_1(I_k)\simeq
H_1(\overline{I}_k)$, hence one may first reduce modulo $J$. Moreover, the
standard short exact sequence
\begin{equation}\label{h1}
{\diagram
0 \rto & \B_1 \rto & \Z_1 \rto & H_1(\overline{I}_k) \rto & 0
\enddiagram }
\end{equation}
implies that $\lambda(H_1(I_k))=\lambda(H_1(\overline{I}_k))=
\lambda(\Z_1)-\lambda(\B_1)$.
Let
\[
\diagram
0 \rto & \Z_1 \rto & (R/J)^{n_k} \rto & I_k/J \rto & 0
\enddiagram
\]
be a minimal presentation of $I_k/J$. Then
\begin{equation}\label{lenz1}
\lambda(\Z_1)=n_k\lambda(R/J)-\lambda(I_k/J)=n_k\lambda(R/I_k)+
(n_k-1)\lambda(R/\m^k),
\end{equation}
as, by $($\ref{11}$)$, $\lambda(I_k/J)=\lambda(R/\m^k)$.
On the other hand $\B_1$ fits in the short exact sequence
\begin{equation}\label{sv-bu}
{ \diagram
0 \rto & \B_1 \rto & (I_k/J)^{n_k} \rto & S_2(I_k/J) \rto & 0,
\enddiagram }
\end{equation}
which is essentially in \cite{SV3} $($see also \cite[proof of Theorem
3.1]{BU}$)$, so that
\[
\lambda(\B_1)=n_k \lambda(I_k/J) - \lambda(S_2(I_k/J))=n_k\lambda(R/\m^k)-\rho.
\]
Finally, $($\ref{h1}$)$, $($\ref{lenz1}$)$, and $($\ref{sv-bu}$)$
give that
\begin{eqnarray*}
\lambda(H_1(I_k)) & = &
n_k\lambda(R/I_k) + (n_k-1)\lambda(R/\m^k) - \left(n_k\lambda(R/\m^k)-
\rho\right) \\
& = & n_k\lambda(R/I_k) + \rho - \lambda(R/\m^k),
\end{eqnarray*}
as in Proposition~\ref{Theorem2.4.8}. }
\end{Remark}

A technical fact is needed in order to prove Theorem~\ref{H_1} below;
even though we only need this result in the case $k=2$, we prove it in
full generality.

\begin{Lemma}\label{lem4}
Let $(R, \m)$ be a Gorenstein local ring of dimension $d\geq 2$
and let $J=(z_1, \ldots, z_d)$ be an ideal generated by a regular
sequence contained in $\m^s$, where $s$ is a positive integer $\geq
2$. For $k=1, \ldots, s$ set $I_k=J \colon \m^k$ and, for $k=3,
\ldots, s$, assume that $\m^k \colon \m = \m^{k-1}$.
Then there exists a non degenerate bilinear form
\[
\diagram
\theta_k \colon\ \m^{k-1}/\m^k \times I_k/I_{k-1} \rto
& R/\m,
\enddiagram
\]
which determines a well-defined, non-singular square matrix whose size is
$n_k=\lambda(\m^{k-1}/\m^k)$.
\end{Lemma}
\begin{pf}
Given $\widetilde{m}=m+\m^k\in \m^{k-1}/\m^k$ and
$\widehat{l}=l+I_{k-1}\in I_k/I_{k-1}$ one has that
$ml\equiv cv\pmod{J}$, as $\m^{k-1} I_k\subseteq I_1$ and $I_1=(J,
v)$. Thus a bilinear form $\theta_k$ can be defined, up to an element of
$\m$, by
\begin{equation}\label{bubba}
\theta_k(\widetilde{m}, \widehat{l})=c.
\end{equation}
It is routine to verify that $\theta_k$ is a well-defined bilinear form.
Let us prove, instead, the non-degeneracy of $\theta_k$. To say that
$\theta_k(\widetilde{m}, \widehat{l}) \equiv 0\pmod{\m}$ for all
$\widehat{l}\in I_k/I_{k-1}$ means
that $ml\in J$ for all $l\in I_k$. But $\m^k=J\colon
I_k$, so that $m\in \m^k$. On the other hand,
$\theta_k(\widetilde{m}, \widehat{l}) \equiv 0\pmod{\m}$ for all
$\widetilde{m} \in \m^{k-1}/\m^k$ implies that
$ml\in J$ for all $m\in \m^{k-1}$; hence $l\in I_{k-1}$, since
$I_{k-1}=J\colon \m^{k-1}$, so that $l \equiv 0\pmod{I_{k-1}}$.

Both $\m^{k-1}/\m^k$ and $I_k/I_{k-1}$ are $R/\m$-vector spaces of
dimension $n_k$; hence, $\theta$ can be written in terms of fixed
bases of those two vector spaces. Let $\{ \widetilde{x}_1, \ldots,
\widetilde{x}_{n_k} \}$ and $\{ \widehat{b}_1, \ldots,
\widehat{b}_{n_k} \}$ be bases of $\m^{k-1}/\m^k$ and
$I_k/I_{k-1}$ respectively. With respect to those two bases let
$\widetilde{m}=(a_1, \ldots a_{n_k})$ and $\widehat{l}=(d_1, \ldots
d_{n_k})$. Finally, letting $x_ib_j \equiv c_{ij}v\pmod{J}$, one has
\begin{eqnarray*}
ml
& \equiv & \left(\sum_{i=1}^{n_k} a_ix_i \right)  \left(
\sum_{j=1}^{n_k} d_jb_j \right) = \sum_{i=1}^{n_k} a_i
\sum_{j=1}^{n_k} d_jx_ib_j  \qquad \bmod J\\
& \equiv & \left(\sum_{i=1}^{n_k} a_i \sum^{n_k}_{j=1}
d_jc_{ij} \right) s.
\end{eqnarray*}
Thus, the $c$ that appears in $($\ref{bubba}$)$ is $\sum_i
a_i\sum_j d_jc_{ij}$; in matrix form, the former expression can be
written as $($note that the $c_{ij}$'s are seen modulo $\m)$
\[
\left(\mbox{$\,a_1\ \cdots \ a_{n_k} \,$}
\right)
\left( \begin{array}{ccc}
c_{11} & \cdots & c_{1 n_k} \\
\vdots &  & \vdots \\
c_{n_k 1} & \cdots & c_{n_k n_k}
\end{array}
\right)
\left(\begin{array}{c}
d_1 \\
\vdots \\
d_{n_k}
\end{array}
\right).
\]
Since $\theta_k$ is a non-degenerate quadratic form, $(c_{ij})$ is a non
singular matrix. \qed
\end{pf}

\begin{Theorem}\label{H_1}
Let $(R, \m)$ be a Gorenstein local ring of dimension $d\geq 2$
and let $J=(z_1, \ldots, z_d)$ be an ideal generated by a regular
sequence contained in $\m^s$, where $s$ is a positive integer $\geq
2$.
For $k=1, 2$ set $I_k=J \colon \m^k$ and
assume that $z_1, \ldots, z_d$ are among the minimal generators of
$I_2$. Then the length of the first Koszul homology module $I_2$ is
\[
\lambda(H_1(I_2))=n\lambda(R/I_2)+ {n+1 \choose 2}- \lambda(R/\m^2),
\]
where $n=n_2=\lambda(\m/\m^2)$.
\end{Theorem}
\begin{pf}
According to Remark~\ref{remh_1}, it will be enough to show that
\begin{equation}\label{len-b1}
\lambda(\B_1)=n\lambda(R/\m^2)-{n+1\choose 2}=n(1+n)-{n+1\choose
2}={n+1\choose 2},
\end{equation}
as $\lambda(R/\m^2)=\lambda(R/\m)+\lambda(\m/\m^2)=1+n$. Note that
$\B_1$ is the submodule of
$\overline{R}\,{}^n$ generated by the $n$-uples $\overline{u}_{ij}$
whose non-zero components are in the $i$-th and in the $j$-th
positions and are given by $\overline{b}_j$ and $-\overline{b}_i$
respectively, i.e.,
\[
\overline{u}_{ij}=(\ldots, \overline{b}_j, \ldots, -\overline{b}_i,
\ldots);
\]
hence, $\mu(\B_1)={n\choose 2}$.
In order to show $($\ref{len-b1}$)$, we use the short exact
sequence
\[
\diagram
0 \rto & \overline{\m}\B_1 \rto & \B_1 \rto &
\B_1/\overline{\m}\B_1 \rto & 0.
\enddiagram
\]
Observe that both $\B_1/\overline{\m}\B_1$ and $\overline{\m}\B_1$
are $R/\m$-vector spaces and that $($\ref{len-b1}$)$ follows provided
one shows that
\[
\lambda(\B_1/\overline{\m}\B_1)={n\choose 2} \qquad
\mbox{and} \qquad \lambda(\overline{\m}\B_1)=n.
\]
However $\lambda(\B_1/\overline{\m}\B_1)=\mu(\B_1)={n\choose
2}$, thus it remains to show that $\lambda(\overline{\m}\B_1)=n$.
As in the proof of Lemma~\ref{lem4}, one can consider
$\overline{\m}=(\overline{x}_1, \ldots, \overline{x}_n)$  and
write $\overline{x_ib}_j=\overline{c_{ij}v}$. Thus a
typical element of $\overline{\m}\B_1$ is of the form
\begin{eqnarray*}
\overline{mu}_{ij} & = & \left(
\sum_{k}\overline{a_kx}_k\right)\overline{u}_{ij}=\sum_k
\overline{a_kx_ku}_{ij} = \sum_k \overline{a}_k(\ldots,
\overline{x_kb}_j, \ldots, -\overline{x_kb}_i, \ldots)   \\
& = & \sum_k \overline{a_kv}(\ldots, \overline{c}_{kj},
\ldots, -\overline{c}_{ki}, \ldots).
\end{eqnarray*}
So $\overline{\m}\B_1$ is a subspace of a $n$-dimensional
$R/\m$-vector space and it is generated by the following
$n{n\choose 2}$ vectors
\[
\overline{v}(\ldots, \overline{c}_{kj}, \ldots, -\overline{c}_{ki},
\ldots),
\]
where $i, j, k=1, \ldots n$. One can identify
$\overline{\m}\B_1$ as a subspace of $(R/\m)^n$ by
identifying the previous vectors with the $n$-uples $($defined modulo
$R/\m)$
\[
(\ldots, c_{kj}, \ldots, -c_{ki}, \ldots).
\]
Now, if the subspace spanned by those vectors was a proper subspace
then it would be contained in some hyperplane, i.e., those vectors would all
satisfy an equation of the form
\begin{equation}\label{hyperplane}
\alpha_1X_1+\alpha_2X_2\cdots+\alpha_nX_n=0,
\end{equation}
where $\alpha=(\alpha_1, \ldots, \alpha_n)$ is a non-trivial $n$-upla. Fix
$i, j$ and substitute the vector $(\ldots, c_{kj}, \ldots,
-c_{ki}, \ldots)$ in $($\ref{hyperplane}$)$ for each $X_k$. This leads
to the equations
\[
\begin{array}{l}
\alpha_1 c_{1j}+\alpha_2 c_{2j}+ \cdots
+\alpha_n c_{nj} =0, \\
\alpha_1 c_{1i}+\alpha_2 c_{2i}+ \cdots
+\alpha_n c_{ni} =0.
\end{array}
\]
Letting $i$ and $j$ vary between $1$ and $n$ yields the
following homogeneous system of linear equations
\[
\left\{
\begin{array}{ccccccccc}
\alpha_1 c_{11}&+&\alpha_2 c_{21}&+&\cdots&
+&\alpha_n c_{n1} &=&0, \\
\vdots && \vdots &&  && \vdots && \vdots\, \\
\alpha_1 c_{1n}&+&\alpha_2 c_{2n}&+&\cdots&
+&\alpha_n c_{nn} &=&0.
\end{array}
\right.
\]
But $\alpha\not=0$ contradicts the non-singularity of the matrix
$(c_{ij})$ proved in Lemma~\ref{lem4}. Thus,
$\lambda(\overline{\m}\B_1)=n$. \qed
\end{pf}

\begin{ack}
The authors wish to acknowledge Wolmer V. Vasconcelos, for useful
discussions they had during the writing of this paper, and Bernd Ulrich,
for helpful suggestions which improved the exposition of Section 3.
\end{ack}

\end{document}